\documentclass[a4paper,11pt,reqno]{amsart}

\usepackage[margin=20mm,includehead,includefoot]{geometry}
\usepackage[utf8]{inputenc}
\usepackage[T1]{fontenc}
\usepackage{lmodern,microtype}
\usepackage{upref}
\usepackage{mathtools,amssymb,amsthm,amsmath}
\usepackage[foot]{amsaddr}
\usepackage{graphicx}
\usepackage{hyperref}
\usepackage{wrapfig}
\usepackage[margin=1cm]{caption}
\usepackage{bookmark}

\graphicspath{{./graphics/}}

\newtheorem{theorem}{Theorem}

\makeatletter
\let\old@setaddresses\@setaddresses
\def\@setaddresses{\bigskip{\parindent 0pt\let\scshape\relax\let\ttfamily\relax\old@setaddresses}}
\makeatother


\hypersetup{
  pdftitle={A book proof of the middle levels theorem},
  pdfauthor={Torsten M\"utze}
}

\title{A book proof of the middle levels theorem}

\author{Torsten M\"utze}
\address[Torsten M\"utze]{Department of Computer Science, University of Warwick, United Kingdom \& Department of Theoretical Computer Science and Mathematical Logic, Charles University, Prague, Czech Republic}
\email{torsten.mutze@warwick.ac.uk}

\thanks{This work was supported by Czech Science Foundation grant GA~22-15272S}

\begin{document}
\begin{abstract}
We give a short constructive proof for the existence of a Hamilton cycle in the subgraph of the $(2n+1)$-dimensional hypercube induced by all vertices with exactly~$n$ or $n+1$ many 1s.
\end{abstract}

\maketitle

The \emph{$n$-dimensional hypercube $Q_n$} is the graph that has as vertices all bitstrings of length~$n$, and an edge between any two bitstrings that differ in a single bit.
The \emph{weight} of a vertex~$x$ of~$Q_n$ is the number of 1s in~$x$.
The \emph{$k$th level} of~$Q_n$ is the set of vertices with weight~$k$.

\begin{theorem}
\label{thm:middle}
For all~$n\geq 1$, the subgraph of~$Q_{2n+1}$ induced by levels~$n$ and~$n+1$ has a Hamilton cycle.
\end{theorem}

Theorem~\ref{thm:middle} solves the well-known \emph{middle levels conjecture}, and it was first proved in~\cite{MR3483129} (see this paper for a history of the problem). 
A shorter proof was presented in~\cite{MR3819051} (12 pages).
Here, we present a proof from `the book'.

\begin{wrapfigure}{r}{0.42\textwidth}
\centering
\includegraphics[page=1]{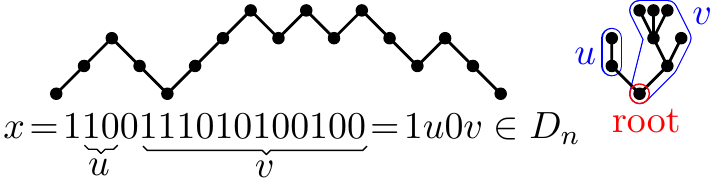}
\caption{A Dyck word (left) and the corresponding ordered rooted tree (right).}
\label{fig:tree}
\end{wrapfigure}

% we don't use the standard proof environment, as it doesn't like wrapfigures
\vspace{2mm}
\textit{Proof.}
We write $D_n$ for all Dyck words of length~$2n$, i.e., bitstrings of length~$2n$ with weight~$n$ in which every prefix contains at least as many 1s as 0s.
We also define $D:=\bigcup_{n\geq 0} D_n$.
Any $x\in D_n$ can be decomposed uniquely as $x=1u0v$ with $u,v\in D$.
Furthermore, Dyck words of length~$2n$ can be identified by ordered rooted trees with $n$ edges as follows; see Figure~\ref{fig:tree}:
Traverse the tree with depth-first search and write a 1-bit for every step away from the root and a 0-bit for every step towards the root.
For any bitstring~$x$, we write $\sigma^s(x)$ for the cyclic right rotation of~$x$ by $s$ steps.
We write $A_n$ and~$B_n$ for the vertices of~$Q_{2n+1}$ in level~$n$ or~$n+1$, respectively, and we define $M_n:=Q_{2n+1}[A_n\cup B_n]$.
For any $x\in D_n$, $b\in\{0,1\}$ and $s\in\{0,\ldots,2n\}$ we define $\langle x,b,s\rangle:=\sigma^s(xb)$.
Note that we have $A_n=\{\langle x,0,s\rangle \mid x\in D_n\wedge 0\leq s\leq 2n\}$ and $B_n=\{\langle x,1,s\rangle\mid x\in D_n\wedge 0\leq s\leq 2n\}$.
Thus, we think of every vertex of~$M_n$ as a triple~$\langle x,b,s\rangle$, i.e., an ordered rooted tree~$x$ with $n$ edges referred to as the \emph{nut}, a bit~$b\in\{0,1\}$, and an integer~$s\in\{0,\ldots,2n\}$ referred to as the \emph{shift}.

\begin{wrapfigure}{r}{0.33\textwidth}
\centering
\includegraphics[page=2]{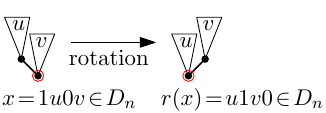}
\caption{Tree rotation.}
\label{fig:rot}
\end{wrapfigure}
The first step is to construct a cycle factor in the graph~$M_n$.
For this we define a mapping $f:A_n\cup B_n\rightarrow A_n\cup B_n$ as follows.
Given an ordered rooted tree $x=1u0v\in D_n$ with $u,v\in D$, a \emph{tree rotation} yields the tree~$r(x):=u1v0\in D_n$; see Figure~\ref{fig:rot}.
We define $f(\langle x,0,s\rangle):=\langle r(x),1,s+1\rangle$ and $f(\langle x,1,s\rangle):=\langle x,0,s\rangle$.
It is easy to see that $f$ is a bijection.
Indeed, the inverse mapping is $f^{-1}(\langle x,0,s\rangle)=\langle x,1,s\rangle$ and $f^{-1}(\langle x,1,s\rangle)=\langle r^{-1}(x),0,s-1\rangle$.
Furthermore, $f$ changes only a single bit.
To see this observe that for $x=1u0v$ with $u,v\in D$ the bitstrings $\langle x,0,s\rangle=\sigma^s(1u0v0)$ and $f(\langle x,0,s\rangle)=\langle r(x),1,s+1\rangle=\sigma^{s+1}(u1v01)=\sigma^s(1u1v0)$ differ only in the bit between the substrings~$u$ and~$v$.
We also note that $f^2(\langle x,0,s\rangle)=\langle r(x),0,s+1\rangle\neq \langle x,0,s\rangle$.
Consequently, for any vertex~$y$ of~$M_n$, the sequence $C(y):=\big(y,f(y),f^2(y),\ldots\big)$ is a cycle, and $F_n:=\{C(y)\mid y\in A_n\cup B_n\}$ is a cycle factor in~$M_n$.

As $f^2(\langle x,0,s\rangle)=\langle r(x),0,s+1\rangle$, moving two steps forward along a cycle of~$F_n$ applies a tree rotation to the nut, and increases the shift by~$+1$.
As the ordered rooted tree~$x\in D_n$ has $n$ edges, we have $x=r^{2n}(x)$.
Consequently, the minimum integer $t>0$ such that $x=r^t(x)$ must divide~$2n$.
It follows that $\gcd(t,2n+1)=1$, hence all shifts of the nut~$x$ are contained in the cycle~$C(\langle x,0,0\rangle)$, i.e., $\langle x,0,s\rangle \in C(\langle x,0,0\rangle)$ for all $s\in \{0,\ldots,2n\}$.
Therefore, the cycles of~$F_n$ are in bijection with equivalence classes of ordered rooted trees with $n$ edges under tree rotation, also known as \emph{plane trees}.
In particular, the number of cycles of~$F_n$ is the number of plane trees with $n$ edges (OEIS~A002995).

\begin{wrapfigure}{r}{0.28\textwidth}
\centering
\includegraphics[page=3]{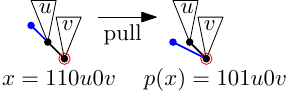}
\caption{Pull operation.}
\label{fig:pull}
\end{wrapfigure}
\begin{wrapfigure}{r}{0.42\textwidth}
\centering
\vspace{-4mm}
\includegraphics[page=4]{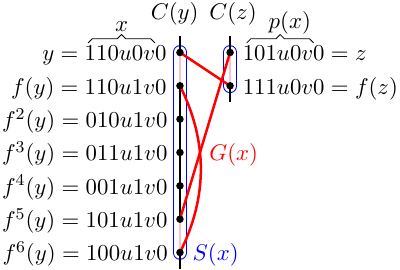}
\caption{Gluing 6-cycle $G(x)$.}
\label{fig:gluing}
\end{wrapfigure}

The second step is to glue the cycles of the factor~$F_n$ to a single Hamilton cycle.
We call an ordered rooted tree~$x\in D_n$ \emph{pullable} if $x=110u0v$ for $u,v\in D$, and we define $p(x):=101u0v\in D_n$.
We refer to $p(x)$ as the tree obtained from~$x$ by a \emph{pull} operation.
In words, the leftmost leaf of~$x$ is in distance~2 from the root, and the edge leading to this leaf is removed and reattached as the new leftmost child of the root in~$p(x)$; see Figure~\ref{fig:pull}.
For any pullable tree~$x=110u0v\in D_n$ with $u,v\in D$, we define $y:=\langle x,0,0\rangle=x0$ and $z:=\langle p(x),0,0\rangle=p(x)0$, and we consider the 6-cycle $G(x):=(y,f(y),f^6(y),f^5(y),z,f(z))=(110u0v0,110u1v0,100u1v0,101u1v0,101u0v0,111u0v0)$, which has the edges $(y,f(y))$ and~$(f^6(y),f^5(y))$ in common with the cycle~$C(y)$, and the edge~$(z,f(z))$ in common with the cycle~$C(z)$; see Figure~\ref{fig:gluing}.
Consequently, if~$C(y)$ and~$C(z)$ are two distinct cycles, then the symmetric difference between the edge sets of~$C(y)$, $C(z)$ and~$G(x)$ is a single cycle on the same set of vertices, i.e., $G(x)$ glues the cycles~$C(y)$ and~$C(z)$ together.

\noindent
We define $S(x):=\{f^i(y)\mid i=0,\ldots,6\}\cup\{z,f(z)\}$, and we claim that for any two pullable trees~$x\neq x'$, we have $S(x)\cap S(x')=\emptyset$, i.e., the cycles~$C(x)$ and~$C(x')$ are (vertex-)disjoint.
To see this, consider the shifts of the vertices in~$S(x)$ and~$S(x')$, which are $0,1,1,2,2,3,3,0,1$.
It follows that if $S(x)\cap S(x')\neq\emptyset$, then we have $x=x'$, $p(x)=x'$, or $x=p(x')$.
These cases are ruled out by the assumption $x\neq x'$, the fact that $p(x)=10\cdots$ and $x'=11\cdots$ differ in the second bit, and that $x=11\cdots$ and $p(x')=10\cdots$ differ in the second bit, respectively.

To complete the proof, it remains to show that the cycles of the factor~$F_n$ can be glued to a single cycle via gluing cycles~$G(x)$ for a suitable set of pullable trees~$x\in D_n$.
As argued before, none of the gluing operations interfere with each other.
Using the interpretation of the cycles of~$F_n$ as equivalence classes of ordered rooted trees under tree rotation, it suffices to prove that every cycle can be glued to the cycle that corresponds to the star with $n$ edges.
As each gluing cycle corresponds to a pull operation, this amounts to proving that any ordered rooted tree~$x\in D_n$ can be transformed to the star~$(10)^n$ via a sequence of tree rotations and/or pulls.

\begin{wrapfigure}{r}{0.37\textwidth}
\centering
\includegraphics[page=5]{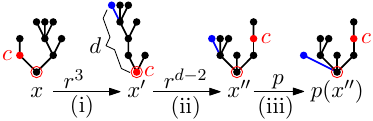}
\caption{Illustration of steps~(i)--(iii) that make a tree more star-like.}
\label{fig:connect}
\end{wrapfigure}
Indeed, this is achieved as follows:
We fix a vertex~$c$ of~$x$ to become the center of the star (this vertex never changes), and we repeatedly perform the following three steps; see Figure~\ref{fig:connect}:
(i) rotate~$x$ to a tree~$x'$ such that~$c$ is root and the leftmost leaf of~$x'$ is in distance~$d>1$ from~$c$; (ii) apply $d-2$ rotations to $x'$ to obtain a tree~$x''$ whose leftmost leaf has distance~2 from the root; (iii) perform a pull.
As step~(iii) decreases the sum of distances of all vertices from~$c$, we reach the star after finitely many steps.

This completes the proof of the theorem.
\qed

\section*{Acknowledgements}

Arturo Merino suggested the triple notation~$\langle x,b,s\rangle$, which allowed further streamlining of the proof.

\bibliographystyle{alpha}
\bibliography{../refs}

\begin{thebibliography}{GMN18}

\bibitem[GMN18]{MR3819051}
P.~Gregor, T.~M\"{u}tze, and J.~Nummenpalo.
\newblock A short proof of the middle levels theorem.
\newblock {\em Discrete Anal.}, Paper No. 8:12~pp., 2018.

\bibitem[M{\"u}t16]{MR3483129}
T.~M{\"u}tze.
\newblock Proof of the middle levels conjecture.
\newblock {\em Proc. Lond. Math. Soc.}, 112(4):677--713, 2016.

\end{thebibliography}

\end{document}